\documentclass[11pt]{article}
\usepackage{amsmath,amssymb,amsthm,amsfonts}
\usepackage{cite}
\usepackage[hidelinks]{hyperref}

\allowdisplaybreaks[4]
\newtheorem{thm}{Theorem}[section]

\newtheorem{lem}[thm]{Lemma}
\newtheorem{prop}[thm]{Proposition}

\theoremstyle{definition}

\theoremstyle{remark}

\numberwithin{equation}{section}

\makeatletter
\renewenvironment{proof}[1][\proofname]{\par
  \pushQED{\qed}%
  \normalfont\topsep6\p@\@plus6\p@\relax
  \trivlist
  \item[\hskip\labelsep\normalfont\bfseries #1\@addpunct{.}]\ignorespaces
}{%
  \popQED\endtrivlist\@endpefalse
}
\makeatother

\newcommand{\R}{\mathbb{R}}
\newcommand{\norm}[1]{\left\lVert #1\right\rVert}
\newcommand{\abs}[1]{\left\lvert #1\right\rvert}
\DeclareMathOperator{\dist}{dist}
\DeclareMathOperator{\supp}{supp}
\DeclareMathOperator{\diver}{div}

\begin{document}

\title{Uniqueness and Nondegeneracy of Positive Bounded Solutions for Semilinear Dirichlet Problems in Exterior-Star Domains}
\author{ Wei Ji\thanks{School of Mathematics and Statistics, and  Key Laboratory of Nonlinear Analysis $\&$ 
Applications (Ministry of Education), Central China Normal University, P.O. Box 71010, Wuhan 430079, P. R. China. 
Email: \texttt{jiwei2020@amss.ac.cn}. }
,\, Shuang Wu\thanks{School of Mathematics and Statistics,  and Hubei Key Laboratory of Mathematical Sciences, 
Central China Normal University, Wuhan 430079, P. R. China. Email: \texttt{swu@mails.ccnu.edu.cn}.}
}
\date{}
\maketitle

\begin{abstract}
We prove the uniqueness and nondegeneracy of positive bounded solutions
for the Dirichlet problem
\begin{equation*}
 -\Delta u=f(u)\,\ \ \text{in}\ \ \Omega, \ \ u=0\,\ \ \text{on}\ \ \partial\Omega,
\end{equation*}
where $f$ is a positive reaction and $\Omega\subsetneq\mathbb R^d$
is a uniformly $C^{2,\gamma}$ domain with a star-shaped complement. Our uniqueness result gives a positive answer to Open Question~1 of [H.~Berestycki and C.~Graham, ARMA, 2026].
\end{abstract}

\medskip

\section{Introduction}

In this paper, we study the positive bounded solutions of the following
Dirichlet problem
\begin{equation}\label{mainequation}
 \begin{cases}
  -\Delta u=f(u)\,\ \ \text{in}\ \ \Omega,\\
  u=0\,\ \ \text{on}\ \ \partial\Omega,
 \end{cases}
\end{equation}
where $\Omega\subsetneq\R^d$ is a nonempty connected open set and $d\geq1$. 
We assume that $\Omega$ is uniformly $C^{2,\gamma}$
for some $\gamma\in(0,1]$ in the sense of
\cite[Definition~A.1]{BG2025}. The nonlinearity satisfies
\begin{equation}\label{eq-reaction}
 \begin{gathered}
f\in C^{1,\gamma}([0,+\infty)),\ \ f(0)=f(1)=0,\ \ f'(0)>0,\ \ f'(1)<0,\\
f(s)>0\ \ \text{for}\ \ 0<s<1,\ \ f(s)<0\ \ \text{for}\ \ s>1.
\end{gathered}
\end{equation}
Following Berestycki and Graham \cite{BG2026},
we call $f$ satisfying \eqref{eq-reaction}
a positive reaction. A positive bounded solution of \eqref{mainequation} is a
function
\begin{equation*}
 u\in C^2(\Omega)\cap C(\overline\Omega)\cap L^\infty(\Omega),\ \ u>0
 \,\ \ \text{in}\ \ \Omega
\end{equation*}
satisfying \eqref{mainequation} pointwise.
When $\Omega$ is unbounded, we impose no additional condition
on $u(x)$ as $\abs{x}\to+\infty$.

Problem \eqref{mainequation} is associated with the following
reaction--diffusion equation with zero Dirichlet boundary
conditions:
\begin{equation}
 \label{eq-parabolic}
 \partial_t U=\Delta U+f(U)
 \,\ \ \text{in}\ \ (0,+\infty)\times\Omega,
 \ \ U=0
 \,\ \ \text{on}\ \ (0,+\infty)\times\partial\Omega.
\end{equation}
In models of population growth and spatial propagation, 
$U$ represents the population density
and $f(U)$ describes local population growth, see 
\cite{Fisher1937,KPP1991,AronsonWeinberger1975,
AronsonWeinberger1978,Skellam1951,CantrellCosner2003}.
For stationary solutions $U(t,x)=u(x)$ of \eqref{eq-parabolic}
where $u=0$ on $\partial\Omega$,
we have
$\partial_tU=0$, so that the problem \eqref{eq-parabolic}
reduces to \eqref{mainequation}.
 
When $\Omega=\R^d$, the constant $u\equiv1$ is the unique
positive bounded solution of $-\Delta u=f(u)$ under
\eqref{eq-reaction}, see \cite[Section~3]{BG2026}.
When $\Omega$ is smooth and bounded, $f$ satisfies
\eqref{eq-reaction}, and $s\mapsto f(s)/s$ is strictly
decreasing on $(0,+\infty)$, Brezis and Oswald
\cite{BO1986} proved that \eqref{mainequation} has at most
one positive solution.
If $f$ satisfies \eqref{eq-reaction} and the strong-KPP
condition that $s\mapsto f(s)/s$ is strictly decreasing on
$(0,1]$, Berestycki and Graham \cite{BG2025} proved that
\eqref{mainequation} has at most one positive bounded solution
even when $\Omega$ is unbounded, provided that $f'(0)$ lies
outside the closure of the principal limit spectrum of the
Dirichlet Laplacian on $\Omega$.
If $f$ satisfies \eqref{eq-reaction} and $\Omega$ is the epigraph
of a globally Lipschitz function, Berestycki, Caffarelli and
Nirenberg \cite{BCN1997} proved that \eqref{mainequation} has
a unique positive bounded solution.

Following \cite{BG2026}, we call $\Omega$
exterior-star if its complement $K:=\R^d\setminus\Omega$ is
nonempty and star-shaped. More precisely, $K$ is star-shaped
about $x_*\in K$ if
\begin{equation*}
 x_*+t(y-x_*)\in K
 \ \ \text{for any}\ \ y\in K \ \ \text{and}\ \ t\in[0,1].
\end{equation*}
The domain $\Omega$ is strongly exterior-star if $K$ is
star-shaped about each point of a nonempty open subset of $K$.
Even a smooth compact set can have only one star center
\cite[Section~2.1]{BG2026}.
Berestycki and Graham \cite[Theorem~2.1]{BG2026} proved that
\eqref{mainequation} has a unique positive bounded solution
if $f$ satisfies \eqref{eq-reaction}, $\Omega$ is strongly
exterior-star, and there exist $R>0$ and a convex set
$C\subset\R^d$ such that
\begin{equation}
 \label{eq-BG-convexity}
 K\setminus B_R(0)=C\setminus B_R(0).
\end{equation} 
They also asked whether the uniqueness of positive bounded
solutions of \eqref{mainequation} still holds when $K$ is
only assumed to be star-shaped
\cite[Section~8.1, Open Question~1]{BG2026}:

\emph{If $f$ is positive, does \eqref{mainequation} have exactly
one positive bounded solution on all exterior-star domains?}

Recall that  ``$f$ is positive'' means that $f$ satisfies
\eqref{eq-reaction}.
We give a positive answer to this question in the following uniqueness
theorem.

\begin{thm}\label{thm-main}
Suppose $d\geq1$, $\gamma\in(0,1]$, $f$ satisfies \eqref{eq-reaction}, and
$\Omega\subsetneq\R^d$ is a nonempty connected uniformly $C^{2,\gamma}$
exterior-star domain. Then \eqref{mainequation} has exactly one positive
bounded solution.
\end{thm}

Theorem~\ref{thm-main} completely answers the open question above. 
Moreover, we point out that in the proof of Theorem~\ref{thm-main}, the
strong star-shapedness of $K$ can be weakened to star-shapedness,
and the convexity condition \eqref{eq-BG-convexity} can be removed.
In the proof of Theorem~\ref{thm-main}, we also show that
for any star center $x_*$ of $K$, the unique positive bounded
solution $u$ of \eqref{mainequation} is strictly increasing
in the radial coordinate centered at $x_*$, in the sense that
\begin{equation}\label{eq-all-centers-monotonicity}
 (x-x_*)\cdot\nabla u(x)>0
 \,\ \ \text{in}\ \ \Omega.
\end{equation}
The positive reaction assumption \eqref{eq-reaction} is essential for uniqueness
in Theorem~\ref{thm-main}. Indeed, when $d\geq2$ and
$f(s)=s(s-\theta)(1-s)$ with $0<\theta<1/2$,
Berestycki, Graham and Wei
\cite{BGW2025} proved that
\eqref{mainequation} has infinitely many positive bounded
solutions even when
$\Omega=\R^d\setminus\overline{B_1(0)}$.

We prove the uniqueness by comparing local energies.  To this end, we first construct the smallest and largest positive bounded solutions  $u_-$ and $u_+$ and prove their strict radial monotonicity in Lemma~\ref{prop-extremal}.  We then adapt the calibration for monotone families of solutions in \cite[Section~4]{AAC2001} 
to the dilations of $u_\pm$, using the formulation for subsolutions and supersolutions in \cite{CEF2024}. This gives that both solutions locally minimize the energy $E_1$ defined in \eqref{eq-local-energy}; see Lemma~\ref{prop-minimality}. The local minimality of $u_+$ next implies that $u_+(x/\kappa)$ minimizes $E_{\kappa^{-2}}$ in $\kappa\Omega$ for every $\kappa>1$. Since $\kappa\Omega\subseteq\Omega$ and $1-\kappa^{-2}>0$, Proposition~\ref{prop-minimizer-order} gives that $u_+(x/\kappa)\leq u_-(x)$ in $\kappa\Omega$. The proof of this comparison uses energy differences on bounded sets and a cutoff argument (cf. \cite{BR2015}). Finally, letting $\kappa\to1^+$ gives that $u_+\leq u_-$ in $\Omega$. Together with $u_-\leq u_+$, this proves that $u_-=u_+$ in $\Omega$.

% The complete writing requirements (updated 30 October 2026) were read
% before preparing this text. The current A1--A30 and G1 were followed.

We also prove the nondegeneracy of the unique positive bounded solution
for \eqref{mainequation}. A positive bounded solution $u$ of
\eqref{mainequation}
is called nondegenerate if the linearized Dirichlet problem
\begin{equation}\label{eq-linearized-problem}
 \begin{cases}
  -\Delta h-f'(u)h=0\,\ \ \text{in}\ \ \Omega,\\
  h=0\,\ \ \text{on}\ \ \partial\Omega
 \end{cases}
\end{equation}
admits only the zero solution in
$C^2(\Omega)\cap C(\overline\Omega)\cap L^\infty(\Omega)$.

\begin{thm}\label{thm-nondegeneracy}
Suppose $d\geq1$, $\gamma\in(0,1]$, $f$ satisfies \eqref{eq-reaction}, and
$\Omega\subsetneq\R^d$ is a nonempty connected uniformly $C^{2,\gamma}$
exterior-star domain. Then the unique positive bounded solution of \eqref{mainequation}
is nondegenerate.
\end{thm}

Without the assumption that $\R^d\setminus\Omega$ is star-shaped,
a unique positive bounded solution of \eqref{mainequation}
can be degenerate even though $f$ satisfies \eqref{eq-reaction}. 
Indeed, Berestycki and Graham (cf. \cite[Proposition~A.5]{BG2026})
constructed a positive reaction $f$ and an interval $\Omega=(0,L)$
for which \eqref{mainequation} has a unique positive bounded
solution $u$ and \eqref{eq-linearized-problem} admits a nonzero
bounded solution $h$.
Under the assumptions of Theorem~\ref{thm-main},
Theorem~\ref{thm-nondegeneracy} excludes such degeneracy.

The key of proving Theorem~\ref{thm-nondegeneracy} is to
show that any bounded solution $h$ of \eqref{eq-linearized-problem}
satisfies
\begin{equation}\label{eq-intro-nd-iteration}
 \int_{\Omega\cap B_R(x_*)}h^2\,dx
 \leq \frac{C}{R}\int_{\Omega\cap B_{2R}(x_*)}h^2\,dx
 \ \ \text{for any}\ \ R\geq1,
\end{equation}
where $x_*$ is a fixed star center of $K$ and $C>0$ is
independent of $R$ and $h$. We obtain this estimate by the
strict radial monotonicity in \eqref{eq-all-centers-monotonicity},
interior and boundary estimates,
and a cutoff argument. Iterating \eqref{eq-intro-nd-iteration}
$d+1$ times and using $h\in L^\infty(\Omega)$, we then obtain that
\begin{equation*}
 \int_{\Omega\cap B_R(x_*)}h^2\,dx
 \leq \frac{C_1}{R}\norm{h}_{L^\infty(\Omega)}^2
 \ \ \text{for any}\ \ R\geq1,
\end{equation*}
where $C_1>0$ is independent of $R$ and $h$.
Letting $R\to+\infty$, we conclude that $h\equiv0$ in $\Omega$.

The paper is organized as follows. In Section~\ref{sec-extremal},
we prove Theorem~\ref{thm-main} by comparing the smallest and
largest positive bounded solutions of \eqref{mainequation}.
In Section~\ref{sec-nondegeneracy}, we establish the proof of
Theorem~\ref{thm-nondegeneracy}.

\section{Extremal solutions and energy comparison}
\label{sec-extremal}

We construct the extremal positive bounded solutions of \eqref{mainequation},
establish their local energy minimality and compare the minimizers to prove
Theorem~\ref{thm-main}.

Let $B_R(x)$ be the open ball of radius $R$ centered at $x$ and set
$B_R:=B_R(0)$. We write $v\in H^1_{\mathrm{loc}}(\overline G)$ if $v\in
H^1(G\cap B_R)$ for every $R>0$. Sobolev boundary values are traces. On the
smooth domains below, zero-trace functions in
$H^1_{\mathrm{loc}}(\overline G)$ extend by zero to
$H^1_{\mathrm{loc}}(\R^d)$.

We first establish the following comparison on unbounded domains.

\begin{lem}\label{lem-linear-comparison}
Suppose $G\subsetneq\R^d$ is open with locally $C^2$ boundary,  and $q\in
H^1_{\mathrm{loc}}(\overline G)\cap L^\infty(G)$ satisfies that for some
constants $c>0$,
\begin{equation*}
 (-\Delta+c)q\leq0\,\ \ \text{in}\ \ G,\ \ q\leq0\,\ \ \text{on}\ \ \partial G
\end{equation*}
in the weak and trace senses. Then $q\leq0$ almost everywhere in $G$.
\end{lem}

\begin{proof}
For $0<\sigma<\sqrt c$, set $\Psi(x):=\sum_{i=1}^d\cosh(\sigma x_i)$ on
$\R^d$. Since $\Psi(x)\to+\infty$ as $\abs{x}\to+\infty$ and $q$ is bounded
with nonpositive trace, $q_\varepsilon:=(q-\varepsilon\Psi)^+\in H^1_0(G)$ has
bounded support for every $\varepsilon>0$. Since $(-\Delta+c)(q
-\varepsilon\Psi)\leq-\varepsilon(c-\sigma^2)\Psi$ in $G$, we have
\begin{equation*}
  0\leq\int_G\bigl(\abs{\nabla q_\varepsilon}^2+c q_\varepsilon^2\bigr)
  dx=\int_G\bigl[\nabla (q-\varepsilon\Psi)\cdot\nabla
  q_\varepsilon+c(q-\varepsilon\Psi)q_\varepsilon\bigr]dx\leq-\varepsilon(c
  -\sigma^2)\int_G\Psi q_\varepsilon dx\leq0.
\end{equation*}
Letting $\varepsilon\to0^+$ then proves Lemma~\ref{lem-linear-comparison}.
\end{proof}

Inspired by \cite[Lemma~2.7]{BG2026}, we next establish the following
lemma.

\begin{lem}\label{lem-deep-interior}
Suppose $f$ satisfies \eqref{eq-reaction}, $G\subsetneq\R^d$ is a nonempty
connected uniformly $C^{2,\gamma}$ domain, and $u$ is a positive bounded
solution of \eqref{mainequation} with $G$ in place of $\Omega$. Then
\begin{equation}\label{eq-range-solution}
 0<u<1\,\ \ \text{in}\ \ G.
\end{equation}
For each $\theta\in(0,1)$, there exists $R_\theta>0$, depending only on
$d,f,\theta$, such that
\begin{equation}\label{eq-deep-interior}
 u(x)\geq\theta\ \ \text{for}\ \ x\in G\ \ \text{with}\ \ \dist(x,G^c)\geq
 R_\theta.
\end{equation}
\end{lem}

\begin{proof}
The boundary regularity (cf. \cite{GT2001}) gives that $u\in
C^1_{\mathrm{loc}}(\overline G)\cap H^1_{\mathrm{loc}}(\overline G)$. Suppose
$S:=\sup_G u>1$. This then implies from \eqref{eq-reaction} that
$c:=\min_{1\leq s\leq S}[-f(s)/(s-1)]>0$ with continuous extension in view of
\eqref{eq-reaction}. Denote $q:=(u-1)^+$. For $\varepsilon>0$, choose a smooth
convex function $\Phi_\varepsilon$ satisfying $\Phi_\varepsilon\equiv0$ on
$(-\infty,1]$, $0\leq\Phi_\varepsilon'\leq1$ on $\R$ and
$\Phi_\varepsilon'\equiv1$ on $[1+\varepsilon,\infty)$. We then obtain from
\eqref{mainequation} that for $0\leq\varphi\in C^\infty_c(G)$,
\begin{equation*}
-\int_G\Phi_\varepsilon(u)\Delta\varphi dx=
\int_G\bigl[\Phi_\varepsilon'(u)f(u)-\Phi_\varepsilon''(u)\abs{\nabla
u}^2\bigr]\varphi dx\leq\int_G\Phi_\varepsilon'(u)f(u)\varphi dx.
\end{equation*}
Since $\Phi_\varepsilon(u)\to q$ and $\Phi_\varepsilon'(u)f(u)\to f(1+q)$
pointwise in $G$ as $\varepsilon\to0^+$, by dominated convergence we obtain
the weak inequality
\begin{equation*}
 (-\Delta+c)q\leq f(1+q)+cq\leq0\,\ \ \text{in}\ \ G,\ \ q\equiv0
 \,\ \ \text{on}\ \ \partial G.
\end{equation*}
This then implies from Lemma~\ref{lem-linear-comparison} that $q\equiv0$ in
$G$, contradicting $S>1$. Since we have
\begin{equation*}
 -\Delta(1-u)+\frac{f(u)}{1-u}(1-u)=0\,\ \ \text{in}\ \ G,\ \ 1-u\equiv1
 \,\ \ \text{on}\ \ \partial G,
\end{equation*}
where the coefficient at $u=1$ is $-f'(1)$, the strong maximum principle gives
that $u<1$ in $G$,
which proves \eqref{eq-range-solution}.

Let $\lambda_1(B_r)$ be the first Dirichlet eigenvalue of $-\Delta$ on $B_r$
and $\phi>0$ be the corresponding  eigenfunction satisfying
$\phi(0)=\max_{B_r}\phi=1$. Since $f'(0)>0$,  for $\theta\in(0,1)$ there
exists $r>0$ such that $\lambda_1(B_r)<\min_{0\leq s\leq\theta}f(s)/s$ with
continuous extension in view of \eqref{eq-reaction}. If $\dist(x,G^c)\geq2r$,
set $\zeta(y):=\theta\phi(y-x)$ on $\overline{B_r(x)}\subset G$. We then have
\begin{equation*}
 -\Delta(a\zeta)=\lambda_1(B_r)a\zeta<f(a\zeta)\,\ \ \text{in}\ \ B_r(x),
 \ \ \forall 0<a\leq1.
\end{equation*}
Since $u>0$ on $\overline{B_r(x)}$, we have $a_*:=\max\{a\in[0,1]:a\zeta\leq
u\,\ \ \text{on}\ \ \overline{B_r(x)}\}>0$. Suppose $a_*<1$. Since
$\zeta\equiv0$ on $\partial B_r(x)$, there exists $y\in {B_r(x)}$ such that
$(u-a_*\zeta)(y)=0$, and hence $(-\Delta(u-a_*\zeta))(y)\leq0$. This then
implies that
\begin{equation*}
 0\geq-\Delta(u-a_*\zeta)(y)=f(a_*\zeta(y))-\lambda_1(B_r)a_*\zeta(y)>0.
\end{equation*}
Hence $u(x)\geq\zeta(x)=\theta$ for $\dist(x,G^c)\geq R_\theta:=2r$.
This proves \eqref{eq-deep-interior}, which completes the proof of
Lemma~\ref{lem-deep-interior}.
\end{proof}

If $K:=\Omega^c$ is star-shaped about $0$, then we have for $\kappa\geq1$,
\begin{equation}\label{eq-dilation-geometry}
K\subset\kappa K,\ \ \kappa\Omega\subset\Omega,\ \ 
\dist(\kappa x,K)\geq\dist(\kappa x,\kappa K)=\kappa\dist(x,K)
\ \ \text{for}\ \ x\in\Omega.
\end{equation}
This then gives that for $\overline{B_r(x_0)}\subset\Omega$,
\begin{equation}\label{eq-large-balls}
 \overline{B_{\kappa r}(\kappa x_0)}=\kappa\overline{B_r(x_0)}
 \subset\kappa\Omega\subset\Omega,\ \ \forall\kappa\geq1.
\end{equation}

\begin{lem}\label{prop-extremal}
Under the assumptions of Theorem~\ref{thm-main}, suppose $0$ is a star center
of $K$. Then there exist positive bounded solutions $u_-,u_+$ of
\eqref{mainequation} such that for every positive bounded solution $u$ of
\eqref{mainequation},
\begin{equation}\label{eq-extremal-order}
 u_-\leq u\leq u_+\,\ \ \text{in}\ \ \Omega.
\end{equation}
Moreover, we have
\begin{equation}\label{eq-extremal-radial}
 x\cdot\nabla u_\pm>0\,\ \ \text{in}\ \ \Omega.
\end{equation}
\end{lem}

\begin{proof}
In view of \eqref{eq-large-balls}, there exists a ball $B\Subset\Omega$ such
that $\lambda_1(B)<f'(0)$ and  the corresponding Dirichlet eigenfunction
$\phi>0$ satisfies $\max_B\phi=1$. Take $\varepsilon\in(0,1)$ such that
\begin{equation}\label{eq-seed-strict}
 f(s)>\lambda_1(B)s\ \ \text{for}\ \ 0<s\leq\varepsilon.
\end{equation}
Extend $w_0:=\varepsilon\phi$ by zero from $B$. It then follows from
\eqref{eq-seed-strict} that for $0\leq\eta\in C^\infty_c(\Omega)$,
\begin{equation*}
 \int_\Omega\nabla w_0\cdot\nabla\eta dx=\lambda_1(B)\int_B w_0\eta
 dx+\varepsilon\int_{\partial B}\partial_\nu\phi\,\eta dS \leq\int_\Omega
 f(w_0)\eta dx.
\end{equation*}
Let $z\in C^2(\Omega)\cap H^1_{\mathrm{loc}}(\overline\Omega)$ satisfy
\begin{equation*}
 -\Delta z\geq f(z),\ \ 0<z\leq1\,\ \ \text{in}\ \ \Omega,\ \ z\geq0
 \,\ \ \text{on}\ \ \partial\Omega.
\end{equation*}
Since $\min_{\overline B}z>0$ and $-\Delta(aw_0)<f(aw_0)$ in $B$ for
$0<a\leq1$ from \eqref{eq-seed-strict}, similar argument of
Lemma~\ref{lem-deep-interior} gives that
\begin{equation}\label{eq-universal-seed}
 w_0\leq z\,\ \ \text{in}\ \ \Omega.
\end{equation}

Set $M:=1+\max_{0\leq s\leq1}\abs{f'(s)}$, so that $s\mapsto Ms+f(s)$
increases from $0$ to $M$ on $[0,1]$. For measurable $0\leq h\leq1$ in
$\Omega$, let $Th$ be the bounded weak solution of
\begin{equation}\label{eq-iteration-map}
 (-\Delta+M)Th=Mh+f(h)\,\ \ \text{in}\ \ \Omega,\ \ Th\equiv0
 \,\ \ \text{on}\ \ \partial\Omega.
\end{equation}
On each nonempty $\Omega\cap B_j$, by the Lax--Milgram theorem,
\eqref{eq-iteration-map} admits a solution $v_j\in H^1_0(\Omega\cap B_j)$.
Since $0\leq Mh+f(h)\leq M$, we have $0\leq v_j\leq1$. Extending $v_j$ by
zero, we then obtain from \eqref{eq-iteration-map} that for $\eta\in
C^\infty_c(\R^d)$,
\begin{equation}\label{eq-resolvent-local-estimate}
 \int_\Omega\eta^2\abs{\nabla v_j}^2 dx\leq4\int_\Omega\abs{\nabla\eta}^2
 dx+2M\int_\Omega\eta^2 dx.
\end{equation}
We can extract a subsequence such that for every $R>0$ and $\eta\in
C^\infty_c(\R^d)$,
\begin{equation*}
 v_j\rightharpoonup Th\,\ \ \text{in}\ \ H^1(B_R),\ \ \eta
 v_j\rightharpoonup\eta Th\,\ \ \text{in}\ \ H^1_0(\Omega)\ \ \text{as}\ \ j
 \to\infty,
\end{equation*}
where  $0\leq Th\leq1$ solves \eqref{eq-iteration-map} and  satisfies
\eqref{eq-resolvent-local-estimate} with $Th$ in place of $v_j$. We have  for
$h_1\leq h_2$, 
\begin{equation*}
 (-\Delta+M)(Th_1-Th_2)=Mh_1+f(h_1)-Mh_2-f(h_2)\leq0\,\ \ \text{in}\ \ \Omega.
\end{equation*}
Lemma~\ref{lem-linear-comparison} gives that $Th_1\leq Th_2$ almost everywhere
in $\Omega$.

Let $q\in H^1_{\mathrm{loc}}(\overline\Omega)$ be continuous  and satisfy
\begin{equation*}
 -\Delta q\leq f(q),\ \ 0\leq q\leq1\,\ \ \text{in}\ \ \Omega,\ \ q\equiv0
 \,\ \ \text{on}\ \ \partial\Omega.
\end{equation*}
It then follows from \eqref{eq-iteration-map} that
\begin{equation*}
 \begin{aligned}
  (-\Delta+M)(Tz-z)=&f(z)+\Delta z\leq0\,\ \ \text{in}\ \ \Omega,\ \ Tz-z\leq0
  \,\ \ \text{on}\ \ \partial\Omega,\\
  (-\Delta+M)(q-Tq)=&-\Delta q-f(q)\leq0\,\ \ \text{in}\ \ \Omega,\ \ q-Tq
  \equiv0\,\ \ \text{on}\ \ \partial\Omega.
 \end{aligned}
\end{equation*}
By Lemma~\ref{lem-linear-comparison}, $Tz\leq z$ and $q\leq Tq$ almost
everywhere in $\Omega$. Taking $q=w_0$ and $z\equiv1$, we then have $w_0\leq
Tw_0$ and $T1\leq1$. Define $u_n:=T^n w_0$ and $v_n:=T^n1$ for $n\geq0$ (cf.
\cite[Sections~6 and~9]{Amann1976}). This thus follows that
\begin{equation*}
 w_0\leq u_n\leq u_{n+1}\leq v_{n+1}\leq v_n\leq1
 \,\ \ \text{a.e.\ in}\ \ \Omega\ \ \text{for}\ \ n\geq0.
\end{equation*}
By interior $W^{2,p}$ estimates with $p>d$ and Schauder estimates (cf.
\cite{GT2001}), we obtain from \eqref{eq-iteration-map} and
\eqref{eq-resolvent-local-estimate} that there exist unique  $0\leq u_-,u_+\in
C^2_{\mathrm{loc}}(\Omega)\setminus\{0\}$ such that for $\eta\in
C^\infty_c(\R^d)$,
\begin{equation*}
 \begin{gathered}
  (u_n,v_n)\to(u_-,u_+)\,\ \ \text{in}\ \ C^2_{\mathrm{loc}}(\Omega)^2
  \ \ \text{as}\ \ n\to\infty,\\
  (\eta u_n,\eta v_n)\rightharpoonup(\eta u_-,\eta
  u_+)\,\ \ \text{in}\ \ H^1_0(\Omega)^2\ \ \text{as}\ \ n\to\infty.
 \end{gathered}
\end{equation*}
Using the boundary $W^{2,p}$ estimates (cf. \cite{GT2001}), this then gives
from \eqref{eq-iteration-map} that
\begin{equation*}
 -\Delta u_\pm=f(u_\pm)\,\ \ \text{in}\ \ \Omega,\ \ u_\pm\in
 C(\overline\Omega),\ \ u_\pm\equiv0\,\ \ \text{on}\ \ \partial\Omega.
\end{equation*}
The strong maximum principle and Lemma~\ref{lem-deep-interior} then give that
$0<u_-\leq u_+<1$ in $\Omega$. Since it follows from \eqref{eq-universal-seed}
that
\begin{equation*}
 u_n=T^n w_0\leq T^nz\leq z,\ \ q\leq T^nq\leq T^n1=v_n
 \,\ \ \text{a.e.\ in}\ \ \Omega\ \ \text{for}\ \ n\geq0,
\end{equation*}
we obtain that $u_-\leq z$ and $q\leq u_+$ in $\Omega$. By
Lemma~\ref{lem-deep-interior} and boundary regularity, we can take $z=q=u$ for
any positive bounded solution $u$ of \eqref{mainequation}. This proves
\eqref{eq-extremal-order}.

Fix $\kappa>1$ and set $z_\kappa(x):=u_-(\kappa x)$ in $\Omega$ and
$w_\kappa(x):=u_+(x/\kappa)$ in $\kappa\Omega$, extending $w_\kappa$ by zero
to $\Omega$. This then implies from \eqref{eq-dilation-geometry} that
$z_\kappa,w_\kappa\in H^1_{\mathrm{loc}}(\overline\Omega)$ and $w_\kappa$ is
continuous with zero trace. Since $-\Delta z_\kappa=\kappa^2 f(z_\kappa)\geq
f(z_\kappa)$ in $\Omega$, we obtain that $u_-(x)\leq z_\kappa(x)=u_-(\kappa x)$
for $x\in\Omega$. For $w_\kappa$, we have $-\Delta
w_\kappa=\kappa^{-2}f(w_\kappa)>0$ in $\kappa\Omega$. By the Hopf lemma, this
implies that $\partial_\nu w_\kappa<0$ on $\partial(\kappa\Omega)$ for the
outward normal $\nu$. We thus obtain that for $0\leq\eta\in
C^\infty_c(\Omega)$,
\begin{equation*}
 \int_\Omega\nabla w_\kappa\cdot\nabla\eta dx=\kappa^{-2}
 \int_{\kappa\Omega}f(w_\kappa)\eta dx+
 \int_{\partial(\kappa\Omega)}\partial_\nu w_\kappa\,\eta dS\leq\int_\Omega
 f(w_\kappa)\eta dx.
\end{equation*}
Taking $q=w_\kappa$ in $q\leq u_+$, we derive that $u_+(x)=w_\kappa(\kappa
x)\leq u_+(\kappa x)$ for $x\in\Omega$. Thus we have
\begin{equation*}
 u_\pm(x)\leq u_\pm(\kappa x)\ \ \text{for}\ \ x\in\Omega
 \ \ \text{and}\ \ \kappa\geq1.
\end{equation*}
This then implies that $Z:=x\cdot\nabla w\geq 0$ in $\Omega$ for $w=u_+$ or
$w=u_-$. It follows from \eqref{mainequation} that
\begin{equation}\label{eq-radial-derivative}
 -\Delta Z-f'(w)Z=2f(w)>0\,\ \ \text{in}\ \ \Omega.
\end{equation}
If $Z(y)=0$ for some $y\in\Omega$, then we obtain from
\eqref{eq-radial-derivative} that $0\geq-\Delta Z(y)=2f(w(y))>0$, a
contradiction. This proves \eqref{eq-extremal-radial}, which completes the
proof of Lemma~\ref{prop-extremal}.
\end{proof}

Define the local energy for a constant $a>0$ and a bounded open set
$A\subset\R^d$,
\begin{equation}\label{eq-local-energy}
 E_a(v;A):=\int_A\left(\frac12\abs{\nabla v}^2-aF(v)\right) dx,\ \ v\in
 H^1(A;[0,1]),
\end{equation}
where $F(s):=\int_0^s f(t)\,dt$ for $0\leq s\leq1$. Let $G\subset\R^d$ be open with locally $C^2$ boundary and  $w\in
H^1_{\mathrm{loc}}(\overline G;[0,1])$ have zero trace. We call $w$ an
$E_a$-minimizer in $G$ if for every $R>0$ and every $\psi\in H^1_0(G)$ with
zero extension supported in $B_R$ and $0\leq w+\psi\leq1$ almost everywhere in
$G$, we have
\begin{equation}\label{eq-minimizer-definition}
 E_a(w+\psi;G\cap B_R)\geq E_a(w;G\cap B_R).
\end{equation}

The following lemma gives a sufficient condition for a positive bounded
solution of \eqref{mainequation} to be an $E_1$-minimizer in $\Omega$.

\begin{lem}\label{prop-minimality}
Under the assumptions of Theorem~\ref{thm-main}, let $K$ be star-shaped about
$0$ and $w$ be a positive bounded solution of \eqref{mainequation} satisfying
$x\cdot\nabla w(x)>0$ in $\Omega$. Then $w$ is an $E_1$-minimizer in $\Omega$.
\end{lem}

\begin{proof}
By Lemma~\ref{lem-deep-interior} and local boundary regularity, we have
$0<w<1$ in $\Omega$ and $w\in H^1_{\mathrm{loc}}(\overline\Omega)$ with zero
trace. It follows from \eqref{eq-dilation-geometry} that $I_x:=\{t\in\R:e^t
x\in\Omega\}=(t_x,+\infty)$ with $-\infty\leq t_x<0$ for $x\in\Omega$, where
we use the convention $e^{-\infty}:=0$. Define
\begin{equation}\label{eq-dilation-family}
 w_t(x):=w(e^t x),\ \ \partial_t w_t(x)=e^t x\cdot\nabla w(e^t x)>0,\ \ t\in
 I_x.
\end{equation}
Since $0\in K$, we obtain from \eqref{eq-dilation-geometry} that
\begin{equation*}
 e^t x\to e^{t_x}x\in\partial\Omega\ \ \text{as}\ \ t\searrow t_x,\ \ 
 \dist(e^t x,K)\geq e^t\dist(x,K)\to+\infty\ \ \text{as}\ \ t\to+\infty.
\end{equation*}
Since $w_t(x)$ is continuous and $w(x)=0$ on $\partial\Omega$, we obtain from
Lemma \ref{lem-deep-interior} that
\begin{equation}\label{eq-family-limits}
 w_t(x)\to0\ \ \text{as}\ \ t\searrow t_x,\ \ w_t(x)\to1\ \ \text{as}\ \ 
 t\to+\infty.
\end{equation}
By \eqref{eq-dilation-family} and \eqref{eq-family-limits}, each
$(x,s)\in\Omega\times(0,1)$ determines a unique $\tau(x,s)\in I_x$ satisfying
\begin{equation}\label{eq-inverse-family}
 w(e^{\tau(x,s)}x)=s,\ \ \operatorname{sgn}\tau(x,s)=\operatorname{sgn}(s
 -w(x)).
\end{equation}
Together with \eqref{eq-dilation-family}, we obtain by the implicit function
theorem that $\tau\in C^1(\Omega\times(0,1))$. Define
\begin{equation*}
 P(x,s):=e^{\tau(x,s)}\nabla w(e^{\tau(x,s)}x),\ \ (x,s)\in\Omega\times(0,1).
\end{equation*}
Since it follows from \eqref{eq-inverse-family} that $P(x,w_t(x))=\nabla
w_t(x)$ in $x\in\Omega\cap e^{-t}\Omega$ at fixed $t$, we obtain from
\eqref{mainequation} that for $(x,s)\in\Omega\times(0,1)$,
\begin{equation}\label{eq-field-divergence}
 \diver_xP(x,s)+\partial_sP(x,s)\cdot P(x,s)=\Delta
 w_t(x)|_{t=\tau(x,s)}=-e^{2\tau(x,s)}f(s).
\end{equation}

Let $A\Subset\Omega$ be a smooth bounded open set and take $\psi\in
H^1_0(\Omega)$ with $\supp\psi\Subset A$. Set $v:=w+\psi$ and suppose that
$\delta_0\leq v,w\leq1-\delta_0$ almost everywhere in $A$ for some
$\delta_0\in(0,1/2)$. Set $B(x):=\int_{w(x)}^{v(x)}P(x,s) ds$ in $A$. Then we
have
\begin{equation}\label{eq-field-integration}
 0=\int_A\diver B(x) dx=\int_A\bigl[P(x,v(x))\cdot\nabla v-|\nabla w|^2\bigr]
 dx+\int_A\int_{w(x)}^{v(x)}\diver_xP(x,s) ds dx,
\end{equation}
where we have used $\supp B\Subset A$ in the first equality and
$P(x,w(x))=\nabla w(x)$ in the second equality. Since
$\frac12|P(x,v(x))|^2-\frac12|\nabla w(x)|^2=\int_{w(x)}^{v(x)}P(x,s)
\cdot\partial_sP(x,s) ds$, we deduce from \eqref{eq-field-divergence} and
\eqref{eq-field-integration} that
\begin{align}
E_1(v;A)-E_1(w;A)=&\frac12\int_A|\nabla v-P(x,v(x))|^2
dx+\int_A\int_{w(x)}^{v(x)}\bigl(e^{2\tau(x,s)}-1\bigr)f(s) ds dx\notag\\
\geq&\int_A\int_{\min\{w(x),v(x)\}}^{\max\{w(x),v(x)\}}\big|e^{2\tau(x,s)}
-1\big|f(s) ds dx\geq0\label{eq-calibration-identity}, 
\end{align}
where we have used \eqref{eq-reaction} and \eqref{eq-inverse-family} in the
first inequality.

Fix $R>0$ and let $\psi$ be admissible in \eqref{eq-minimizer-definition}. We
may assume that $\psi\not\equiv0$, since \eqref{eq-minimizer-definition} holds
with equality for $\psi\equiv0$. Set $d_\Omega(x):=\dist(x,\partial\Omega)$
for $x\in\Omega$. For sufficiently small $r\in(0,1)$, finitely many normal
charts cover $\Omega\cap B_R\cap\{d_\Omega<r\}$. In these charts, we have
\begin{equation*}
 x=y+t\nu(y),\ \ d_\Omega(x)=-t,\ \ dx=J(y,t) dS(y) dt,\ \ 0<c\leq J(y,t)\leq
 C,
\end{equation*}
where $y\in\partial\Omega\cap\overline{B_{R+1}}$, $-r<t<0$ and $\nu$ is the
outward unit normal. Define $ g(t):=\psi(y+t\nu(y))\in  H^1(-r,0)$ with
$g(0)=0$ for almost every $y\in\partial\Omega$.
By Hardy's inequality $\int_{-r}^0|g(t)/t|^2 dt\leq4\int_{-r}^0|g'(t)|^2 dt$,
we obtain that
\begin{equation}\label{eq-local-hardy}
\begin{aligned}
\int_{\Omega\cap B_R\cap\{d_\Omega<r\}}\frac{|\psi(x)|^2}{d_\Omega(x)^2}dx
&\leq C\int_{\partial\Omega\cap\overline{B_{R+1}}}
\int_{-r}^0\frac{|g(t)|^2}{t^2}\,dtdS(y)\\
&\leq 4C\int_{\partial\Omega\cap\overline{B_{R+1}}}\int_{-r}^0|\nabla\psi(y
+t\nu(y))\cdot\nu(y)|^2dtdS(y)\\
&\leq C'\int_\Omega|\nabla\psi(x)|^2dx<\infty.
\end{aligned}
\end{equation}
For $0<\delta<r/2$, set $\eta_\delta:=\min\{1,\max\{0,d_\Omega/\delta-1\}\}$
and $v_\delta:=w+(1-\delta)\eta_\delta\psi$ in $\Omega$. Set $v:=w+\psi$.
Since $\eta_\delta\equiv1$ on $\{x\in\Omega:d_\Omega(x)\geq2\delta\}$ and
$|\nabla\eta_\delta|\leq2/d_\Omega$ almost everywhere in $\Omega$, we obtain
from \eqref{eq-local-hardy} that as $\delta\to0^+$,
\begin{equation}\label{eq-boundary-approximation}
\begin{aligned}
\|v_\delta-v\|_{H^1(\Omega)}^2\leq&2\|(\eta_\delta-1)\psi\|_{H^1(\Omega)}^2
+2\delta^2\|\psi\|_{H^1(\Omega)}^2\\
\leq& C\int_{\Omega\cap B_R\cap\{d_\Omega<2\delta\}}\Big(|\psi|^2
+|\nabla\psi|^2+\frac{|\psi|^2}{d_\Omega^2}\Big) dx
+2\delta^2\|\psi\|_{H^1(\Omega)}^2\to0.
\end{aligned}
\end{equation}
Since $\eta_\delta\equiv0$ on $\{x\in\Omega:d_\Omega(x)\leq\delta\}$, we have
$\supp(v_\delta-w)\Subset\Omega\cap B_R$. There then exists a smooth bounded
open set $A_\delta\Subset\Omega\cap B_R$ containing $\supp(v_\delta-w)$. Since
$0\leq(1-\delta)\eta_\delta(x)\leq 1-\delta$ and $v_\delta=[1-(1
-\delta)\eta_\delta]w+(1-\delta)\eta_\delta(x)v$, we have
\begin{equation*}
 0<\delta\min_{\overline{A_\delta}}w\leq v_\delta,w\leq1
 -\delta\min_{\overline{A_\delta}}(1-w)<1\ \ \text{for a.e.}\ \ x\in A_\delta.
\end{equation*}
From \eqref{eq-calibration-identity} and \eqref{eq-boundary-approximation}, we
obtain that
\begin{equation*}
 \begin{aligned}
  0\leq& E_1(v_\delta;A_\delta)-E_1(w;A_\delta)=E_1(v_\delta;\Omega\cap
  B_R)-E_1(w;\Omega\cap B_R)\\
  \to& E_1(v;\Omega\cap B_R)-E_1(w;\Omega\cap B_R)\ \ \text{as}\ \ 
  \delta\to0^+,
 \end{aligned}
\end{equation*}
where we have used the Lipschitz continuity of $F$ on $[0,1]$. This proves
\eqref{eq-minimizer-definition} with $a=1$, which completes the proof of
Lemma~\ref{prop-minimality}.
\end{proof}

We next compare two local minimizers, when one domain contains the other and
the energy on the larger domain has the larger potential coefficient.

\begin{prop}\label{prop-minimizer-order}
Suppose $f$ satisfies \eqref{eq-reaction}, $D\subseteq G\subsetneq\R^d$ are
nonempty open sets with locally $C^2$ boundaries and $0<a<b$. Suppose $w$ is
an $E_a$-minimizer in $D$ and $u$ is an $E_b$-minimizer in $G$, as defined by
\eqref{eq-minimizer-definition}. Extend $w$ by zero to $G$. Then
\begin{equation}\label{eq-minimizer-order}
 w\leq u\,\ \ \text{almost everywhere in}\ \ G.
\end{equation}
\end{prop}

\begin{proof}
Extend $w,u$ by zero to $\R^d$ and set $h:=(w-u)^+\in
H^1_{\mathrm{loc}}(\R^d)$. Define for compactly supported Lipschitz
$0\leq\chi\leq1$ on $\R^d$, 
\begin{equation}\label{eq-exchange-competitors}
 u_\chi:=u+\chi h,\ \ w_\chi:=w-\chi h\,\ \ \text{in}\ \ \R^d.
\end{equation}
Since $0\leq\chi h\leq h\leq1$ and $h\equiv0$ almost everywhere outside $D$,
we have
\begin{equation*}
 \chi h\in H^1_0(D)\subset H^1_0(G),\ \ 0\leq u_\chi,w_\chi\leq1
 \,\ \ \text{a.e.\ in}\ \ \R^d.
\end{equation*}

It follows from \eqref{eq-exchange-competitors} that
\begin{equation*}
 \begin{aligned}
  Q_\chi :=& \frac12\bigl(\abs{\nabla u_\chi}^2-\abs{\nabla u}^2+\abs{\nabla
  w_\chi}^2-\abs{\nabla w}^2\bigr)\\
  =&-\chi(1-\chi)\abs{\nabla h}^2+(2\chi-1)h\nabla
  h\cdot\nabla\chi+h^2\abs{\nabla\chi}^2.
 \end{aligned}
\end{equation*}
This then implies that
\begin{equation}\label{eq-gradient-cost}
 Q_\chi\leq\frac{\abs{\nabla\chi}^2}{4\chi(1-\chi)}h^2
 \,\ \ \text{a.e.\ on}\ \ \{0<\chi<1\},\ \ Q_\chi\equiv0
 \,\ \ \text{a.e.\ on}\ \ \{\chi=0\}\cup\{\chi=1\}.
\end{equation}
Since $F(s)-Ls^2/2$ is concave on $[0,1]$ for $L:=\max_{[0,1]}\abs{f'}$, we have almost everywhere on
$\{h>0\}$
\begin{equation}\label{eq-potential-exchange}
 \begin{aligned}
V_\chi :=&-b\bigl[F(u+\chi h)-F(u)\bigr]-a\bigl[F(w-\chi h)-F(w)\bigr]\\
\leq&-b\chi[F(w)-F(u)]-a\chi[F(u)-F(w)]+\frac{(a+b)L}{2}\chi(1-\chi)h^2\\
=&-(b-a)\chi[F(w)-F(u)]+\frac{(a+b)L}{2}\chi(1-\chi)h^2,
\end{aligned}
\end{equation}
On $\{h=0\}$, we have $V_\chi\equiv0$ almost
everywhere. Using the continuous extension in view of \eqref{eq-reaction}, set
$c_f:=\frac16\min_{0\leq r\leq1}f(r)/[r(1-r)]>0$. It then follows that almost
everywhere on $\{h>0\}$,
\begin{equation}\label{eq-potential-gap}
 F(w)-F(u)\geq 6c_f\int_u^w r(1-r) dr=6c_f\left[hu(1-w)+\frac{h^2}{2}
 -\frac{h^3}{3}\right]\geq c_fh^2.
\end{equation}
We thus obtain from  \eqref{eq-gradient-cost}--\eqref{eq-potential-gap} that
\begin{equation}\label{eq-interior-gain}
 Q_\chi+V_\chi\leq -(b-a)c_fh^2\,\ \ \text{a.e.\ on}\ \ \{\chi=1\}.
\end{equation}

Define for $R>1$, 
\begin{equation*}
 \chi_R(x):=\begin{cases}
  1, & \abs{x}\leq R,\\
  \displaystyle\cos^2\!\left(\frac\pi2(\abs{x}-R)\right), & R<\abs{x}<R+1,\\
  0, & \abs{x}\geq R+1,
 \end{cases}
\end{equation*}
so that $\abs{\nabla\chi_R}^2/[4\chi_R(1-\chi_R)]=\pi^2/4$ and
$\chi_R(1-\chi_R)\leq1/4$ in $B_{R+1}\setminus\overline{B_R}$. By the
minimality of $u,w$, we obtain from
\eqref{eq-gradient-cost}, \eqref{eq-potential-exchange},
\eqref{eq-potential-gap} and \eqref{eq-interior-gain} that
\begin{equation*}
\begin{aligned}
 0\leq&  E_b(u_{\chi_R};G\cap B_{R+2})-E_b(u;G\cap B_{R+2})+
 E_a(w_{\chi_R};D\cap B_{R+2})-E_a(w;D\cap B_{R+2})\\
 =&\int_{\R^d}(Q_{\chi_R}+V_{\chi_R}) dx \leq -(b-a)c_f\int_{B_R}h^2
 dx+C_{a,b,f}\int_{B_{R+1}\setminus B_R}h^2 dx.
 \end{aligned}
\end{equation*}
where $C_{a,b,f}:=\pi^2/4+(a+b)L/8$. This then implies that for
$m(R):=\int_{B_R}h^2 dx$ with $R>1$,
\begin{equation}\label{eq-mass-recurrence}
 m(R+1)\geq \left(1+\frac{(b-a)c_f}{C_{a,b,f}}\right)m(R).
\end{equation}
If $m(R_0)>0$ for some $R_0>1$, we then deduce from \eqref{eq-mass-recurrence}
and $h\leq1$ that
\begin{equation*}
 0<m(R_0)\leq\frac{m(R_0+n)}{\left(1+(b-a)c_f/C_{a,b,f}\right)^n}
 \leq\frac{\abs{B_1}(R_0+n)^d}{\left(1+(b-a)c_f/C_{a,b,f}\right)^n}\to0
 \ \ \text{as}\ \ n\to\infty,
\end{equation*}
where we have used $0\leq h\leq1$ in the last inequality. Thus $h\equiv0$
almost everywhere, which gives \eqref{eq-minimizer-order} and completes the
proof of Proposition~\ref{prop-minimizer-order}.
\end{proof}

\subsection{Proof of Theorem~\ref{thm-main}}
\label{sec-order}
In this subsection, we mainly establish the proof of Theorem~\ref{thm-main}.
\begin{proof}[Proof of Theorem~\ref{thm-main}]
Translate a star center of $\Omega^c$ to the origin.
It follows from Lemma~\ref{prop-extremal} that there exist
positive bounded solutions $u_\pm$ of \eqref{mainequation}
satisfying \eqref{eq-extremal-radial}.
Applying Lemma~\ref{prop-minimality}, we then obtain that
$u_\pm$ are $E_1$-minimizers in $\Omega$.

Fix $\kappa>1$ and set $w_\kappa(x):=u_+(x/\kappa)$ in
$D:=\kappa\Omega\subset\Omega$ from \eqref{eq-dilation-geometry}. For $\psi$
admissible in \eqref{eq-minimizer-definition} for $w_\kappa$ with zero
extension supported in $B_{\kappa R}$, we have
\begin{equation*}
 \psi(\kappa\cdot)\in H^1_0(\Omega),\ \ \supp\bigl(\psi(\kappa\cdot)\bigr)
 \subset B_R,\ \ 0\leq u_++\psi(\kappa\cdot)\leq1
 \,\ \ \text{a.e.\ in}\ \ \Omega.
\end{equation*}
By the minimality of $u_+$, we obtain from \eqref{eq-local-energy} and
\eqref{eq-minimizer-definition} that
\begin{equation*}
 \begin{aligned}
  E_{\kappa^{-2}}(w_\kappa+\psi;D\cap B_{\kappa R}) =& \kappa^{d-2}
  E_1(u_++\psi(\kappa\cdot);\Omega\cap B_R)\\
  \geq& \kappa^{d-2}E_1(u_+;\Omega\cap B_R)= E_{\kappa^{-2}}(w_\kappa;D\cap
  B_{\kappa R}).
 \end{aligned}
\end{equation*}
This gives that $w_\kappa$ is an $E_{\kappa^{-2}}$-minimizer in $D$. Applying
Proposition~\ref{prop-minimizer-order} with $w=w_\kappa$ and $u=u_-$, we
deduce from \eqref{eq-minimizer-order} and continuity that
\begin{equation}\label{eq-final-dilation-order}
 u_+(x/\kappa)\leq u_-(x)\,\ \ \text{in}\ \ \kappa\Omega,\ \ \forall \kappa>1.
\end{equation}
For fixed $x\in\Omega$, we have $x/\kappa\in\Omega$ for $\kappa>1$
sufficiently close to $1$. Letting $\kappa\to1^+$ in
\eqref{eq-final-dilation-order}, we obtain from \eqref{eq-extremal-order} that
every positive bounded solution $u$ of \eqref{mainequation} satisfies
\begin{equation*}
 u_+(x)\leq u_-(x)\leq u(x)\leq u_+(x).
\end{equation*}
This completes the proof of Theorem~\ref{thm-main}.
\end{proof}

\section{Nondegeneracy of the positive bounded solution}
\label{sec-nondegeneracy}

The purpose of this section is to prove
Theorem~\ref{thm-nondegeneracy}.
Let $u$ be the unique positive bounded solution of
\eqref{mainequation} obtained by Theorem~\ref{thm-main}.
Define the quadratic form associated with
\eqref{eq-linearized-problem} by
\begin{equation}\label{eq-nd-quadratic-form}
 \mathcal Q[v]
 :=\int_\Omega\bigl(\abs{\nabla v}^2-f'(u)v^2\bigr)\,dx
 \ \ \forall\,v\in H^1_0(\Omega).
\end{equation}
We first derive the following uniform $L^2(\Omega)$ estimate for all functions
in $H^1_0(\Omega)$ with bounded support.

\begin{lem}\label{lem-nd-coercivity}
Under the assumptions of Theorem~\ref{thm-main},
suppose $0$ is a star center of $\Omega^c$.
Then there exists $C>0$, independent of $R$ and $v$, such that
for any $R\geq1$ and $v\in H^1_0(\Omega)$
satisfying $\supp v\subset\overline{B_R}$,
\begin{equation}\label{eq-nd-coercivity}
 \norm{v}_{L^2(\Omega)}^2\leq CR\mathcal Q[v].
\end{equation}

\end{lem}

\begin{proof}
It follows from Theorem~\ref{thm-main} and
Lemma~\ref{prop-extremal} that $x\cdot\nabla u(x)>0$
in $\Omega$. We may therefore define
\begin{equation}\label{eq-nd-ZW}
 Z(x):=x\cdot\nabla u(x),\ \ W(x):=\frac{2f(u(x))}{Z(x)}
 \,\ \ \text{in}\ \ \Omega.
\end{equation}
Since $0<u<1$ in $\Omega$, we deduce from
\eqref{eq-reaction} and \eqref{eq-radial-derivative} that
\begin{equation*}
 -\Delta Z-f'(u)Z=2f(u)>0,\ \ W>0
 \,\ \ \text{in}\ \ \Omega.
\end{equation*}
Integrating by parts in \eqref{eq-nd-quadratic-form}, we then
deduce from \eqref{eq-radial-derivative} and \eqref{eq-nd-ZW} that
\begin{equation}\label{eq-nd-form-identity}
 \mathcal Q[v]
 =\int_\Omega Z^2\left|\nabla\left(\frac{v}{Z}\right)\right|^2\,dx
  +\int_\Omega Wv^2\,dx, 
 \ \ v\in C^\infty_c(\Omega).
\end{equation}
For any $v\in H^1_0(\Omega)$, choose $v_j\in C^\infty_c(\Omega)$ such that
$v_j\to v$ in $H^1_0(\Omega)$ and $v_j(x)\to v(x)$ for almost every
$x\in\Omega$ as $j\to\infty$.
Since $f'(u)\in L^\infty(\Omega)$, we have
$\mathcal Q[v_j]\to\mathcal Q[v]$ as $j\to\infty$.
Applying Fatou's lemma, we then deduce from
\eqref{eq-nd-form-identity} that
\begin{equation}\label{eq-nd-weighted-form}
 \mathcal Q[v]\geq\int_\Omega Wv^2\,dx\geq0
 \ \ \forall\,v\in H^1_0(\Omega).
\end{equation}

We next derive an $L^2$ estimate for $v$ near $\partial\Omega$.
Define $d_\Omega(x):=\dist(x,\partial\Omega)$ in $\Omega$.
Since $\Omega$ is uniformly $C^{2,\gamma}$, the principal
curvatures of $\partial\Omega$ are uniformly bounded.
Choose $\rho_0>0$ sufficiently small such that the normal coordinates
\begin{equation*}
 x=y+t\nu(y),\ \ y\in\partial\Omega,\ \ -\rho_0<t<0
\end{equation*}
give a one-to-one parametrization of
$\{x\in\Omega:d_\Omega(x)<\rho_0\}$ with
$d_\Omega(y+t\nu(y))=-t$,
where $\nu$ denotes the outward unit normal to $\partial\Omega$.
We thus obtain that for any $y\in\partial\Omega$ and $-\rho_0<t<0$,
\begin{equation*}
 dx=J(y,t)\,dS(y)\,dt,\ \ C_0^{-1}\leq J(y,t)\leq C_0,
\end{equation*}
where $C_0\geq1$ is independent of $y$ and $t$. By  the Cauchy--Schwarz inequality, we then have for any $v\in C^\infty_c(\Omega)$, $y\in\partial\Omega$ and $t\in(-\rho_0,0)$,
\begin{equation}\label{eq-nd-normal-estimate}
 \abs{v(y+t\nu(y))}^2\leq -t\int_t^0\abs{\nabla v(y+s\nu(y))}^2\,ds.
\end{equation}
Integrating \eqref{eq-nd-normal-estimate} over
$y\in\partial\Omega$ and $t\in(-\rho,0)$
and using the density of $C^\infty_c(\Omega)$
in $H^1_0(\Omega)$, we then obtain that
for any $v\in H^1_0(\Omega)$ and $0<\rho<\rho_0$,
\begin{equation}\label{eq-nd-boundary-layer}
 \int_{\{d_\Omega<\rho\}}v^2\,dx
 \leq C_b\rho^2\int_\Omega\abs{\nabla v}^2\,dx,
\end{equation}
where $C_b>0$ is independent of $\rho$ and $v$.

We now derive a lower bound for $W(x)$ away from $\partial\Omega$.
Define $M:=\max_{0\leq s\leq1}\abs{f'(s)}$ and choose
$\rho\in(0,\rho_0)$, independently of $R$ and $v$, such that
\begin{equation}\label{eq-nd-rho-choice}
 C_b\rho^2M\leq\frac12.
\end{equation}
We also define $\Omega_\rho:=\{x\in\Omega:d_\Omega(x)\geq\rho\}$.
Note from Lemma~\ref{lem-deep-interior} that $0<u<1$ in $\Omega$.
It follows from \eqref{mainequation} that
\begin{equation}\label{eq-nd-interior-equations}
 -\Delta u=\frac{f(u)}{u}\,u,\ \ 
 -\Delta(1-u)=-\frac{f(u)}{1-u}(1-u)
 \,\ \ \text{in}\ \ \Omega.
\end{equation}
Since $f(0)=f(1)=0$, we get that
$\abs{f(u)/u},\abs{f(u)/(1-u)}\leq M$ in $\Omega$.
For any $x\in\Omega_\rho$, we have $B_{\rho/2}(x)\Subset\Omega$.
Applying the interior Harnack inequality and gradient estimates to
$u$ and $1-u$ in \eqref{eq-nd-interior-equations} on $B_{\rho/2}(x)$
(cf.~\cite{GT2001}), we then obtain that
\begin{equation}\label{eq-nd-interior-gradient}
 \abs{\nabla u(x)}
 \leq C_\rho\min\{u(x),1-u(x)\}
 \ \ \forall\,x\in\Omega_\rho,
\end{equation}
where $C_\rho>0$ is independent of $x$.
Note from \eqref{eq-reaction} that
\begin{equation*}
 c_*:=\inf_{0<s<1}\frac{f(s)}{\min\{s,1-s\}}>0.
\end{equation*}
Since $Z>0$ and $f(u)\geq c_*\min\{u,1-u\}$ in $\Omega$,
we deduce from \eqref{eq-nd-ZW} and \eqref{eq-nd-interior-gradient} that
\begin{equation}\label{eq-nd-Z-upper}
 0<Z(x)\leq\abs{x}\abs{\nabla u(x)}
 \leq\frac{C_\rho}{c_*}\abs{x}f(u(x)),
 \ \ x\in\Omega_\rho.
\end{equation}
It then follows from \eqref{eq-nd-ZW} and \eqref{eq-nd-Z-upper} that
there exists $c_\rho>0$, independent of $R$ and $v$, such that
\begin{equation}\label{eq-nd-interior-weight}
 W(x)\geq\frac{c_\rho}{1+\abs{x}}
 \ \ \forall\,x\in\Omega_\rho.
\end{equation}

Suppose $R\geq1$ and $v\in H^1_0(\Omega)$ satisfies
$\supp v\subset\overline{B_R}$.
We derive from \eqref{eq-nd-quadratic-form},
\eqref{eq-nd-weighted-form}, \eqref{eq-nd-boundary-layer},
\eqref{eq-nd-rho-choice}
and \eqref{eq-nd-interior-weight} that
\begin{equation}\label{eq-nd-combined-estimate}
 \begin{aligned}
  \norm{v}_{L^2(\Omega)}^2
  &\leq\int_{\Omega_\rho}v^2\,dx
       +C_b\rho^2\int_\Omega\abs{\nabla v}^2\,dx\\
  &\leq\frac{1+R}{c_\rho}\mathcal Q[v]
       +C_b\rho^2\bigl(\mathcal Q[v]
       +M\norm{v}_{L^2(\Omega)}^2\bigr)\\
  &\leq\left(\frac{1+R}{c_\rho}+C_b\rho^2\right)\mathcal Q[v]
       +\frac12\norm{v}_{L^2(\Omega)}^2. 
 \end{aligned}
\end{equation}
Since $R\geq1$, we conclude from \eqref{eq-nd-combined-estimate}
that \eqref{eq-nd-coercivity} holds, where the constant $C>0$
is independent of $R$ and $v$.
This completes the proof of Lemma~\ref{lem-nd-coercivity}.
\end{proof}

We are now ready to prove Theorem~\ref{thm-nondegeneracy}.

\begin{proof}[Proof of Theorem~\ref{thm-nondegeneracy}]
Suppose $h\in C^2(\Omega)\cap C(\overline\Omega)\cap L^\infty(\Omega)$
satisfies \eqref{eq-linearized-problem}.
Translate a star center of $\Omega^c$ to the origin.
Applying the boundary regularity lemma
(cf. \cite[Lemma~6.18]{GT2001}) to
\eqref{eq-linearized-problem}, we then obtain that
$h\in H^1_{\mathrm{loc}}(\overline\Omega)$ with zero trace on
$\partial\Omega$.
Choose $\eta\in C^\infty_c(B_2;[0,1])$ satisfying
$\eta\equiv1$ in $B_1$, and define $\eta_R(x):=\eta(x/R)$ in $\R^d$
for any $R\geq1$, where $\eta$ is extended by zero outside $B_2$.
Then $\eta_Rh,\eta_R^2h\in H^1_0(\Omega)$ and
$\supp(\eta_Rh)\subset\overline{B_{2R}}$.
Moreover, we have
\begin{equation}\label{eq-nd-cutoff-properties}
 \eta_R\equiv1\,\ \ \text{in}\ \ B_R,\ \ 
 \abs{\nabla\eta_R}\leq\frac{\norm{\nabla\eta}_{L^\infty(\R^d)}}{R}
 \,\ \ \text{in}\ \ \R^d.
 \end{equation}
Testing \eqref{eq-linearized-problem} with $\eta_R^2h$, we then obtain that
\begin{equation}\label{eq-nd-cutoff}
 \mathcal Q[\eta_Rh]
 =\int_\Omega h^2\abs{\nabla\eta_R}^2\,dx.
\end{equation}
Define for any $R\geq1$, 
\begin{equation*}
 m(R):=\int_{\Omega\cap B_R}h^2\,dx.
\end{equation*}
Applying Lemma~\ref{lem-nd-coercivity} to $\eta_Rh$ with radius $2R$,
we then derive from \eqref{eq-nd-cutoff-properties} and
\eqref{eq-nd-cutoff} that
\begin{equation}\label{eq-nd-mass-recurrence}
 m(R)\leq\norm{\eta_Rh}_{L^2(\Omega)}^2
 \leq 2CR\mathcal Q[\eta_Rh]
 \leq\frac{C_1}{R}m(2R),\ \ R\geq1,
\end{equation}
where $C_1>0$ is independent of $R$ and $h$.

On the contrary, suppose $h\not\equiv0$ in $\Omega$. Then there exists $R_0\geq1$ such that $m(R_0)>0$.
Since $h\in L^\infty(\Omega)$, we have
\begin{equation}\label{eq-nd-mass-bound}
 m(R)\leq\abs{B_1}R^d\norm{h}_{L^\infty(\Omega)}^2,\ \ R\geq1.
\end{equation}
Iterating \eqref{eq-nd-mass-recurrence} $n$ times for any integer $n\geq1$,
we then obtain from \eqref{eq-nd-mass-bound} that
\begin{equation*}
 \begin{aligned}
0<m(R_0)
&\leq\left(\frac{C_1}{R_0}\right)^n2^{-n(n-1)/2}m(2^nR_0)\\
&\leq\abs{B_1}R_0^d\norm{h}_{L^\infty(\Omega)}^2\left(\frac{2^dC_1}{R_0}\right)^n2^{-n(n-1)/2}\to0\ \ \text{as}\ \ n\to\infty, 
 \end{aligned}
\end{equation*}
which contradicts the fact that $m(R_0)>0$.
This implies that $h\equiv0$ in $\Omega$, which completes the proof of
Theorem~\ref{thm-nondegeneracy}.
\end{proof}
\noindent
\textbf{Data Availability}: No data were used for the research described in the article.

\noindent
\textbf{Conflict of interest}: The authors declare no conflict of interest.

\noindent
\textbf{Declaration of generative AI and AI-assisted technologies in the writing process:} During the preparation of
this work, the authors used ChatGPT-6.1 Sol to assist with calculations, the review of mathematical arguments
and calculations, and language editing, including grammar and spelling corrections. After using this tool, the
authors reviewed and edited the content as needed and take full responsibility for the content of the published article.

\end{document}